\documentclass[a4paper,12pt]{amsart}

\usepackage{amssymb,amsbsy,amsmath,amsfonts,amssymb,amscd}
\usepackage{latexsym}
\usepackage{graphics}
\usepackage{color}
\input xy
\xyoption{all}

\theoremstyle{plain}
\newtheorem{thm}{Theorem}[section]

\newtheorem{lemma}[thm]{Lemma}

\theoremstyle{definition}

\newtheorem{example}[thm]{Example}

\newtheorem{question}[thm]{Question}

\numberwithin{equation}{section}
\newcommand{\sA}{{\mathcal A}}

\newcommand{\sE}{{\mathcal E}}

\newcommand{\sI}{{\mathcal I}}

\newcommand{\sP}{{\mathcal P}}

\newcommand{\sR}{{\mathcal R}}

\newcommand{\PP}{\ensuremath{\mathbb{P}}}

\newcommand{\ZZ}{\ensuremath{\mathbb{Z}}}

\newcommand{\hol}{\ensuremath{\mathcal{O}}}

\newcommand\al{\alpha}

\newcommand\De{\Delta}

\newcommand\de{\delta}
\newcommand\e{\epsilon}

\newcommand{\ra}{\ensuremath{\rightarrow}}

\def\eea{\end{eqnarray*}}
\def\bea{\begin{eqnarray*}}

\newcommand\dual{\mathrel{\raise3pt\hbox{$\underline{\mathrm{\thinspace d
\thinspace}}$}}}
\newcommand\qe{\ifhmode\unskip\nobreak\fi\quad $\Box$}       % box for QED

\def\BOX{\hfill\lower.5\baselineskip\hbox{$\Box$}}
\newtheorem{theo}{Theorem}[section]

\newtheorem{remarkk}[theo]{Remark}
\newenvironment{rem}{\begin{remarkk}\rm}{\end{remarkk}}

\newtheorem{prop}[theo] {Proposition}
\newtheorem{cor}[theo]{Corollary}

\newcommand{\Proof}{{\it Proof. }}
\title [Canonical surfaces]{Canonical surfaces of higher degree}

\author{Fabrizio Catanese}
\address {Lehrstuhl Mathematik VIII\\
Mathematisches Institut der Universit\"at Bayreuth\\
NW II,  Universit\"atsstr. 30\\
95447 Bayreuth}
\email{fabrizio.catanese@uni-bayreuth.de}

\thanks{
The present work took place in the realm of the 
 ERC Advanced grant n. 340258, `TADMICAMT' }

\date{\today}

\begin{document}

\maketitle

{\em  Dedicated to Philippe Ellia on the occasion of his 60th birthday}

\begin{abstract}
We consider a family of surfaces of general type $S$ with $K_S$ ample, having  $K^2_S = 24, p_g (S) = 6, q(S)=0$.
We prove  that for these surfaces  the canonical system is base point free and yields an embedding $\Phi_1 : S \ra \PP^5$.

This result answers  a question posed by G. and M. Kapustka \cite{kk}.

We  discuss some related open problems, concerning also the case  $p_g(S) = 5$, where one requires the canonical map
to be birational onto its image.
\end{abstract}

\tableofcontents

\section*{Introduction}

Among the many questions that one may ask about surfaces of general type the following one has not yet been 
sufficiently considered.

\begin{question}\label{q}
Let $S$ be a smooth surface with ample canonical divisor $K_S$,  assume that $p_g(S) = 6$
and that the canonical map $\phi_1 : S \ra \PP^5$ is a biregular embedding.

 Which values of $K^2_S$ can occur,
in particular which  is the maximal value that $K^2_S$ can reach?
\end{question}

Recall the Castelnuovo inequality, holding if $\phi_1$ is birational (onto its image):
$$  K^2_S \geq 3 p_g(S) - 7,$$ 
and the  Bogomolov-Miyaoka-Yau inequality
$$  K^2_S \leq 9 \chi(S) = 9 + 9 p_g(S) - 9 q(S).$$
By virtue of these inequalities,  under the assumptions of  question \ref{q} one must have:
 $$ 11 \leq K^2_S \leq 63.$$
 
 In the article \cite{AMS} methods of homological algebra were used to construct such surfaces with low degree $11 \leq K^2_S \leq 17$,
 and also   to attempt a classification of them. Recently, M. and G. Kapustka constructed in \cite{kk} such canonical surfaces of degree $K^2_S =18$,
 using the method of bilinkage.
 In a preliminary version of the article they even ventured to ask whether the answer to question \ref{q} would be  $K^2_S \leq 18$.
 
 Our main result consists in  exhibiting   such  canonically embedded surfaces having degree $K^2_S = 24$, and with $q(S) = 0$.
 The family of surfaces was indeed  listed in the article  \cite{singular}, dedicated to applications  of the
 technique of bidouble covers; but it was not a priori clear that their  canonical system would be an embedding
 (this was proven in \cite{modsp} for bidouble covers satisfying  much stronger conditions).
 
 \begin{theo}\label{embedding}
Assume that $S$ is a bidouble cover of the quadric $ Q = \PP^1 \times \PP^1$ branched on  three    curves $D_1, D_2, D_3$
of respective bidegrees $(2,3), (2,3) (4,1)$, which 
are smooth and intersect transversally. Assume   moreover that the 12 intersection points $D_1 \cap D_2$ have pairwise different
images  via the second projection $p_v : Q \ra \PP^1$ (i.e., they have different  coordinates  $(v_0 : v_1)$).

Then $S$ is a simply connected \footnote{This follows from theorem 3.8 of \cite{modsp}.} surface with $K^2_S = 24, p_g (S) = 6, q(S)=0$, whose  canonical map
$ \phi_1 : S \ra \PP^5 $ 
is a biregular  embedding. These surfaces form a (non empty)  irreducible algebraic subset    of dimension $25$ of the moduli space.
\end{theo}

We plan  to try to describe the equations of the above  surfaces $S \subset \PP^5$; since,   by a theorem of Walter \cite{walter},
each $S$ is the Pfaffian  locus of a twisted  antisymmetric map $ \al : \sE (-t)  \ra \sE $ of a vector bundle $\sE$ on $\PP^5$,
a natural question is to describe the bundle $\sE$.

Concerning  question \ref{q}, we should point out that it is often easier to construct algebraic varieties 
as parametric images rather than as zero sets of ideals.  Note that 
surfaces with $p_g(S) = 6, q(S) = 2, K^2_S = 45$ (and with $K_S$ ample since they are ball quotients) were constructed in \cite{BC}.
For these, however, the canonical system has base points.

It would be interesting to see whether there do exist  canonically embedded surfaces with  $K^2_S = 56$ which are  regular surfaces isogenous to a product
(see \cite{isogenous}).

We finish this introduction pointing out that a similar question is wide open also for $p_g(S) = 5$ (while for $p_g(S) = 4$ some work has been done, see for instance
 \cite{singular}).

\begin{question}
Let $S$ be a smooth minimal surface of general type,  assume that $p_g(S) = 5$
and that the canonical map $\phi_1 : S \ra \PP^4$ is birational. What is the maximal value that $K^2_S$ can reach?

\end{question}

Again, Castelnuovo's inequalty gives $ 8 \leq K^2_S $, and Bogomolov-Miyaoka-Yau's inequality  gives $ K^2_S \leq 54$.

The known cases where $\phi_1$ is  an embedding are just for $ K^2_S = 8,9$,
where $S$ is a complete intersection of type $(4,2)$ or $(3,3)$. 

Indeed these are the only cases by virtue of the following (folklore?)  theorem which was stated and proven in \cite{AMS}, propositions 6.1 and 6.2, corollary 6.3.

\begin{theo}\label{4}
Assume that $S$ is the minimal model of a surface of general type with $p_g(S) = 5$, and assume that the canonical map $\phi_1 $ embeds 
$S$ in $\PP^4$.

Then $S$ is a complete intersection with $\hol_S (K_S) = \hol_S(1)$, i.e., $S$  is a complete intersection in $\PP^4$ of type $(2,4 )$ or  $(3,3)$.

Moreover, if  $\phi_1 $ is birational, and $K^2_S = 8,9$, then $\phi_1 $ yields an embedding of the canonical model of $X$ 
as a complete intersection in $\PP^4$ of type $(2,4 )$ or  $(3,3)$.

\end{theo}

Recall that the first  main ingredient of proof  is the well known  Severi's  double point formula
(\cite{severi}, see \cite{commalg} for the proof of some transversality claims made by   Severi).

Write  the double point formula in the form stated in  \cite{hartshorne}, Appendix A, 4.1.3. It  gives, once we set   $d : =  K^2_S $:
$$ 12 \chi(S) = (17-d) d .$$
Since $\chi(S) \geq 1$, we get $ 8 \leq d \leq 16$, and viewing the  double point formula as  an equation among integers,
we see that it  is only solvable for 
$ d=8,9 \Rightarrow \chi(S) = 6  \Rightarrow q(S) = 0$, or for $ d=12, \Rightarrow \chi(S) = 5  \Rightarrow q(S) = 1$.

One sees that  the last case cannot occur, since 
 the Albanese map of $S$, 
$ \al : S \ra A$, has as image   an elliptic curve $A$; if one denotes  by $g$ the genus of the Albanese fibres, the slope inequalities
for fibred surfaces of Horikawa, or Xiao, or Konno 
 (\cite{hor}, \cite{xiao}, \cite{konno}), give
$g=2$: hence   $\phi_1$ cannot be birational.

One could replace in theorem \ref{4} the hypothesis that   $\phi_1 $ is an embedding of 
$S$ in $\PP^4$ by the weaker condition that  $\phi_1 $ yields an embedding of the canonical model of $X$ 
via an extension of  the Severi double formula to the case of  surfaces with rational double points as singularities 
\footnote{for which however we have not yet found a reference.}.

\section{The construction of the family of surfaces}

We consider the family of algebraic surfaces  $S$,  bidouble covers (Galois covers with group $(\ZZ/2)^2$) of the quadric
$ Q : = \PP^1 \times \PP^1$, described in the fourth line of page 106 of \cite{singular}. 

This means that we take three divisors $$D_j = \{ \de_j = 0\},   \de_1, \de_2 \in H^0 (\hol_Q (2,3)), \ \de_3 \in H^0 (\hol_Q (4,1)),$$
we consider divisor classes $L_1=  L_2, L_3$ with  $$ \hol_Q (L_1) = \hol_Q (L_2) =  \hol_Q (3,2) , \  \hol_Q (L_3) =  \hol_Q (2,3)$$
and we define  the surface $S$ as  $$ Spec (\hol_Q \oplus w_1 \hol_Q (-L_1) \oplus w_2 \hol_Q (- L_2)  \oplus w_3 \hol_Q (-L_3) ),$$
where the ring structure is given by (here $\{i,j,k\} $ is any permutation of $\{1,2,3\} $)
$$ w_i^2 = \de_j \de_k, \ w_i w_j = w_k \de_k.  $$

We refer to \cite{modsp} and \cite{singular} for the basics on bidouble covers, which show that on the surface $S$ the sections $w_i$ can be written
as a product of square roots of the sections $\de_1, \de_2, \de_3$:
$$  w_i = y_j y_k, \ y_i^2 = \de_i.$$ 

As is in \cite{singular}, page 102, define
$$ N : = 2 K_Q + \sum_1^3 D_i = 2 K_Q + \sum_1^3 L_i ,$$
obtaining a description of the canonical ring of $S$, which is a smooth surface  with ample canonical divisor if the three curves $D_1, D_2, D_3$
are smooth and intersect transversally.
The canonical ring is defined as usual by:

$$ \sR: = \oplus_{m=0}^{\infty} \sR_m, \  \sR_m : =  H^0 (\hol_S (m K_S)),$$
and we have (ibidem)
$$   \sR_{2m+1} = y_1 y_2 y_3 H^0 (\hol_Q (K_Q + m N) ) \oplus (\oplus_{i=1}^3 y_i H^0 (\hol_Q (K_Q + m N + L_i))) $$ 
$$   \sR_{2m} = H^0 (\hol_Q ( m N) ) \oplus (\oplus_{i=1}^3 w_i H^0 (\hol_Q ( m N - L_i)) ).$$

In particular, since in this case $\hol_Q ( N) = \hol_Q (4,3) $, 
$$  H^0 (\hol_S ( K_S)) = \oplus_{i=1}^3 y_i H^0 (\hol_Q (K_Q  + L_i) ) = < y_1 u_0, y_1 u_1, y_2 u_0, y_2 u_1, y_3 v_0, y_3 v_1>,$$
where $u_0, u_1$ is a basis of  $ H^0 ( \hol_Q (1,0) )$, $v_0, v_1$ is a basis of  $ H^0 ( \hol_Q (0,1) )$.

In particular, $p_g(S) = 6$; moreover,  $q(S) = 0$ since $ H^1 ( \hol_Q  ) = H^1 ( \hol_Q (- L_i) ) = 0 , \forall i=1,2,3,$
as follows easily from the K\"unneth formula.

Finally, if $\pi : S \ra Q$ is the bidouble cover, then $2 K_S = \pi^* ( N)$, hence $ K_S$ is ample and $K_S^2 = N^2 = 24$.

\section{Proof that the canonical map is an embedding}

\begin{theo}\label{embedding}
Assume that the three  branch curves  $D_1, D_2, D_3$
are smooth and intersect transversally, and moreover that the 12 intersection points $D_1 \cap D_2$ have pairwise different
images  via the second projection $p_v : Q \ra \PP^1$ (i.e., they have different  coordinates  $(v_0 : v_1)$).

Then the canonical map
$ \phi_1 : S \ra \PP^5, $ $ \phi_1 (x) =$ $$ = (y_1 (x) u_0(x) :  y_1 (x) u_1(x) :  y_2 (x)  u_0(x)  :  y_2 (x) u_1(x) :  y_3 (x)  v_0 (x) : y_3 (x) v_1(x) ),$$
is a biregular  embedding, i.e. we have an isomorphism $$\phi_1 : S \ra \Sigma : = \phi_1 ( S).$$
\end{theo}

\Proof
Let $R_i \subset S, R_i : = \{ y_i =0\}$. The curve $R_i$ maps to $D_i$ with degree $2$, and $R_1 \cap R_2 \cap R_3 = \emptyset$
since by our assumption $D_1 \cap D_2 \cap D_3 = \emptyset$.

{\bf Claim 1:} the canonical system is base-point free.

 In fact, for each point $x$ there is $u_i$, $i \in \{0,1\}$, such that $u_i (x) \neq 0$,
and similarly  $v_j $, $j \in \{0,1\}$, such that $v_j  (x) \neq 0$. Hence $x$ is a base point if $y_1=y_2=y_3=0$ at $x$,
 contradicting $R_1 \cap R_2 \cap R_3 = \emptyset$. 

{\bf Claim 2:} $ \phi_1 $ is a local embedding.

2.1) At the points $x \in R_h \cap R_k$
the two sections $y_h,  y_k$ yield local coordinates, hence our assertion.

2.2) For the  points $x \in S \setminus R$, where $R = R_1 \cup R_2 \cup R_3 $ is the ramification divisor,
let us consider the rational map $F_u : \Sigma \ra \PP^1$, 
induced by the linear projections $$(x_1 : x_2 : x_3 : x_4 : x_5 : x_6 ) \dasharrow (x_1 : x_2 )$$
$$(x_1 : x_2 : x_3 : x_4 : x_5 : x_6 ) \dasharrow (x_3 : x_4 ),$$
 respectively the rational map $F_v : \Sigma \ra \PP^1$, 
 induced by the linear projection $$(x_1 : x_2 : x_3 : x_4 : x_5 : x_6 ) \dasharrow (x_5 : x_6 ).$$
 
 We have that $ (u_0(x) :  u_1(x)) = p_u \circ \pi = F_u \circ \phi_1,  (v_0(x) :  v_1(x)) = p_v \circ \pi = F_v \circ \phi_1$.
 Moreover, $F_u \circ \phi_1$ is a morphism outside the finite set $R_1 \cap R_2$,
 $F_v \circ \phi_1$ is a morphism outside $R_3$. 
 
 Hence
 $\pi = (F_u \circ \phi_1) \times (F_v  \circ \phi_1) : S \setminus R$ is
 a morphism of maximal rank, so $\phi_1$ is a local embedding outside of $R$.
 
 2.3) Set for convenience $ u : = p_u \circ \pi , v : = p_v \circ \pi $.
 
 At the points of $R_1 \setminus(  R_2 \cup R_3)$, then $y_1$ , $u$, $v$ give local coordinates, so we are done;
 similarly for the points of $R_2 \setminus(  R_1 \cup R_3)$.
 
 At the points of $R_3 \setminus(  R_1 \cup R_2)$, we need to show that  $y_3$ and $u$ give local coordinates.
Here, we make the remark that $D_3$ is a divisor of bidegree $(4,1)$, hence $ p_u : D_3 \ra \PP^1$
is an isomorphism. Hence $\de_3, u$ give local coordinates at the points of $D_3$, and
we infer that $y_3$ , $u$, give local coordinates at the points of $R_3 \setminus(  R_1 \cup R_2)$.

{\bf Claim 3:} $ \phi_1 $ is injective. 

Assume that $\phi_1 (x) = \phi_1 (x')$. We observe preliminarily that this condition implies that,
if $x \in R_i$, then also $x' \in R_i$, since for instance $ x \in R_1 \Leftrightarrow x_1 (\phi_1 (x)) = x_2 (\phi_1 (x)) = 0$.

3.1) Assume that $x, x' \in S \setminus R$. 

Then $\pi (x) = \pi (x')$, and since the elements of Galois group $G = (\ZZ/2)^2 = \{ \pm 1\}^3 / \{ \pm 1\}$
are determined by $ \e \in  \{ \pm 1\}^3$, and each of them acts on
$\phi_1 (x) = (x_1 : x_2 : x_3 : x_4 : x_5 : x_6 )$ by 
$$  (x_1 : x_2 : x_3 : x_4 : x_5 : x_6 )  \mapsto    ( \e_1 x_1 : \e_1 x_2 : \e_2 x_3 : \e_2 x_4 : \e_3 x_5 : \e_3 x_6 ),$$
the action of $G$ on $\phi_1 (x) $ has an orbit of cardinality $4$, which is precisely the cardinality of $\pi^{-1}(\pi (x))$.

Hence $\phi_1 (x) = \phi_1 (x') \Rightarrow x = x'.$

3.2) Let $x, x' \in R_1 \setminus(  R_2 \cup R_3)$. Again $\pi (x) = \pi (x')$, and we see that the orbit of $\phi_1 (x)$ under $G$
is the set of points $$  ( 0 : 0 : \e_2 x_3 : \e_2 x_4 : \e_3 x_5 : \e_3 x_6 ),$$ which has cardinality $2$; this is precisely the cardinality of 
$ \pi^{-1} (\pi (x))$, so again we are done.

The case $x, x' \in R_2 \setminus(  R_1 \cup R_3)$ is completely analogous. 

3.3)  Let $x, x' \in R_3 \setminus(  R_1 \cup R_2)$. Again, the orbit of $\phi_1 (x)$ under $G$
has cardinality $2$, so we are done if we show that $\pi (x) = \pi (x')$.  On the set $R_3$, however, $\pi$ is not a morphism,
so we argue differently.

We use instead  that $ u (x) = u(x')$, and that $ p_u : D_3 \ra \PP^1$
is an isomorphism to conclude that  $\pi (x) = \pi (x')$.

3.4) If $x, x' \in R_3 \cap R_2$ (for  $x, x' \in R_3 \cap R_1$ the argument is entirely similar),
we proceed as follows.

We obtain again $ u (x) = u(x')$, and since  $ p_u : D_3 \ra \PP^1$
is an isomorphism, we get that  $\pi (x) = \pi (x')$. However, the restriction of $\pi$ to  $R_3 \cap R_2$  is 
bijective, hence $x = x'$.

3.5) Assume finally that $x, x' \in R_1 \cap R_2$ and use again that the restriction of $\pi$ to  $R_1 \cap R_2$ $\pi$ is 
bijective. It suffices therefore to show that $\pi (x) = \pi (x')$.

The condition $\phi_1 (x) = \phi_1 (x') =    ( 0 : 0 : 0 : 0 : y_3 v_0 : y_3 v_1 ),$ since $y_3$ does not vanish on $R_1 \cap R_2$,
implies that $ v (x) = v(x')$, so the points $\pi (x) ,  \pi (x') \in D_1 \cap D_2$ have the same $v$ coordinate.

Hence, by our assumption, $\pi (x) =  \pi (x')$, exactly as desired.

\qed 

\begin{prop}
The hypotheses of theorem \ref{embedding} define a non empty family of dimension $ 11 + 11 + 9= 31$.
\end{prop}

\Proof
$D_1, D_2$ vary in a linear system of dimension $ 3\times 4 -1= 11$, $D_3$ varies in a linear system of dimension $ 5 \times 2 -1= 9$.

The condition that the divisors intersect transversally is a consequence of the fact  that each of the three linear systems embeds $Q$ in a projective space.

The final condition amounts to the following: write
$$  \de_j = u_0^2 A_j (v) + u_0 u_1 B_j (v) + u_1^2 C_j (v) = 0, \ j=1,2.$$ 

We can view the coefficients of $\de_1, \de_2$,  six degree three homogeneous polynomials in $ v = (v_0 : v_1)$,
 as giving a map $\psi$ of $\PP^1$ inside  the $\PP^5$ with coordinates
$$ (A_1 : B_1 : C_1 : A_2 : B_2 : C_2)$$
parametrizing pairs of homogeneous polynomials of degree 2 in $ u = (u_0 : u_1)$,
$$  \sP_j = u_0^2 A_j  + u_0 u_1 B_j  + u_1^2 C_j  = 0, \ j=1,2.$$

Let $\De$ be the resultant $$\De : = Res_u (\sP_1, \sP_2) = (A_1 C_2 - A_2 C_1)^2 + (B_2 C_1 - C_2 B_1) (A_1 B_2 - A_2 B_1).$$

The resultant  $\De$  defines a hypersurface of degree 4 in $\PP^5$, which is reduced and irreducible, being a non degenerate quadric 
in the variables $B_1, B_2$.

We can therefore choose the six polynomials in a general way so that the  twisted cubic curve $\psi (\PP^1)$ 
 intersects $\De$ transversally in 12 distinct points. 

We conclude because  the image of  $D_1 \cap D_2$   via the projection $p_v$  is given by the 12 zeros of 
$$  f (v_0 : v_1) : = Res_u (\de_1, \de_2) = \De (A_1 (v) : B_1 (v) : C_1 (v) : A_2 (v) : B_2 (v) : C_2 (v)),$$
and by our general choice these are 12 distinct points.

\qed

\begin{rem}
i) Taking into account the group of automorphisms of $Q= \PP^1 \times \PP^1$, we see that the  above family  gives a locally closed
subset of dimension $25$ inside the moduli space of surfaces of general type. This dimension is more than the expected
dimension $ 10 \chi (S) - 2 K^2_S = 70 - 48 = 22$.

ii) since $H^0 (\hol_Q (D_i- L_i) = 0$ for each $i=1,2,3$, there are no natural deformations (\cite{modsp}, \cite{singular});
it is not clear  that our family yields an irreducible component of the moduli space, since the elementary method of 
cor. 2.20 and theorem 3.8 of \cite{modsp} do not apply.

\end{rem}

\section{The canonical ring as  a module over the symmetric algebra $\sA : = Sym (\sR_1)$}

Let us first of all look at the homomorphism $m_2 : Sym^2 (\sR_1) \ra \sR_2$, keeping track of the eigenspace decompositions,
and using the notation $ H^0 (\hol_Q (a,b)) = : V(a,b)$. We have then 

$$  \sR_1 =  0 \oplus y_1 V(1,0) \oplus y_2 V(1,0)  \oplus y_3 V(0,1) ,$$

$$   \sR_2 = V(4,3) \oplus y_2 y_3  V(1,1) \oplus y_1 y_3  V(1,1)  \oplus y_1 y_2  V(2,0) .$$

Since $ V(1,0) \otimes V(0,1 ) \cong V(1,1),$ and we have a surjection  $V(1,0) \otimes V(1,0) \ra V(2,0)$,
the non trivial character spaces of  $   \sR_2$ are in the image of $m_2 : Sym^2 (\sR_1) \ra \sR_2$.

Moreover, the kernel of $y_1 V(1,0) \otimes y_2 V(1,0) \ra y_1 y_2 V(2,0)$ is 1-dimensional,
and provides a quadric $\{ q(x) = 0\}$ containing the canonical image $\Sigma$ of $S$. To simplify our notation,
we directly assume that $ S \subset \PP^5$, via the canonical embedding.

Then the quadric that we obtain is:  $ q(x) : = x_2 x_3 - x_1 x_4$.

The trivial character space of $   \sR_2$ is isomorphic to $  V(4,3) $ and contains  the image of the subspace 
$$ W: =  y_1^2 V(2,0) \oplus  y_2^2 V(2,0) \oplus   y_3^2 V(0,2) \subset Sym^2 (\sR_1),$$
which maps onto 
$$  W' : = \de_1 V(2,0) + \de_2 V(2,0) + \de_3 V(0,2).$$

\begin{lemma}
$W \ra W'$ is an isomorphism.
\end{lemma}

\Proof
Assume that there is a kernel: then there are bihomogeneous polynomials $P_1, P_2, P_3$
such that $ P_1 \de_1 + \de_2 P_2 = P_3 \de_3$.

Since $\de_3$ does not vanish at the 12 points where $\de_1=\de_2 = 0$, $P_3= P_3 (v)$ should
vanish on the projections of these 12 points under $p_v$. But this is a contradiction, since $P_3$
has degree 2, while the 12 projected points are distinct. Hence $P_3 \equiv 0$, 
and $\de_1 | P_2$, a contradiction again since $P_2 $ has bidegree $(2,0)$.

\qed 

\begin{cor}
$S$ is  contained in a unique quadric $\{ q(x) = 0\}$, and $m_2 : Sym^2 (\sR_1) \ra \sR_2$ has image of dimension
$ 21-1 = 20$ and codimension $11$.
\end{cor}
\Proof
$\sR_m$ has dimension $ dim ( \sR_m ) = \chi(S) + \frac{1}{2} m (m-1) K^2_S = 7 + 12 m (m-1)$,
which, for $m=2$, is equal to $31$.

\qed

We shall  now choose $z_1, \dots, z_{11} \in \sR_2 $ which, together with $ Im (m_2)$,   generate $\sR_2$.
The elements $z_1, \dots, z_{11}$  induce a basis of the quotient $\sR_2 / Im (m_2)$.

\begin{theo}
Let $\sA : = Sym (\sR_1)$ be the coordinate ring of $\PP^5$, and consider $\sR$ as an $\sA$-module.
Then $1, z_1, \dots, z_{11}$ is a minimal graded system of generators of $\sR$ as an $\sA$-module. 
\end{theo}

\Proof
In view of the previous observations, it suffices to show that these elements generate $\sR$.

Observe preliminarily that $V(a,b) \otimes V(c,d) \ra V(a+c, b+d)$ is always surjective as soon as $a,b,c,d \geq 0$.

For $\sR_3$, let us write: 

$$ \sR_3  = y_1 y_2 y_3 V(2,1) \oplus y_1 V(5,3)   \oplus y_2 V(5,3) \oplus y_3 V(4,4).    $$

The last three eigenspaces are in the image of $ V(4,3) \otimes   \sR_1 \subset \sR_2 \otimes \sR_1 \ra \sR_3$.

Also the first summand is in the image of $y_1 V(1,0) \otimes y_2 y_3 V(1,1)$.

The same argument works for $\sR_{2m+1}$, while for $\sR_{2m+ 2}$ we find  surjections 
$$  H^0(\hol_Q(  N)) \otimes w_i H^0(\hol_Q( m N - L_i) ) \ra  w_i H^0 ( \hol_Q( (m+1) N -  L_i  ) ),$$ 
and 
$$  H^0(\hol_Q(  N)) \otimes  H^0(\hol_Q( m N ) ) \ra   H^0 ( \hol_Q( (m+1) N ) ).$$

Hence the claimed result follows  by induction on $m$. 

\qed

\begin{rem}
1) Since $ dim \  \sR_3 = 79$,  $ dim \  \sA_3 = 56$, we see that the 12 minimal generators of the module admit
$ 56 + 6 \times 11 - 79 = 43$ relations in degree $3$. 

The module $\sR$ is Cohen-Macaulay, hence it has a length 3 minimal graded resolution.

2) In general,  $ dim  \ \sR_m =    7 + 12 m (m-1)$,  
 $$ dim \  \sA_m =   \begin{pmatrix}
m + 5\\
m
\end{pmatrix},$$
hence the image of $ \sA_m \ra \sR_m$ has dimension less than or equal (because $S$ is contained in a quadric)
  to
  $ dim  \sA_m -  dim  \sA_{m-2}$.
  
  Hence the Hartshorne-Rao module  
 $ \oplus _m H^1 (\sI_S (m)$ is non zero in all degrees $m=2,3,4,5,6$. This shows that the bundle $\sE$ contructed by Walter
 using the Horrocks correspondence should be rather interesting.
\end{rem}

\end{document}